\documentclass{amsart} 
\usepackage{amssymb}
\usepackage{amsmath}
\usepackage{euscript}
\begin{document} 
\input epsf.sty

\title[On the equation $P(f)=Q(g)$]{On the equation $P(f)=Q(g),$ where $P,Q$ are polynomials and $f,g$ are entire functions 
}
 
\author{F. Pakovich}
\thanks{Research supported by the ISF, Grant No. 979/05} 
\email{pakovich@math.bgu.ac.il}
\date{} 
\address{Department of Mathematics\\
Ben Gurion University of the Negev \\
P.O.B. 653 Beer Sheva \\ 84105
Israel }

\def\be{\begin{equation}}
\def\ee{\end{equation}}
\def\bs{$\square$ \vskip 0.2cm}
\def\d{{\rm d}} 
\def\D{{\rm D}} 
\def\I{{\rm I}} 
\def\C{{\mathbb C}} 
\def\N{{\mathbb N}} 
\def\P{{\mathbb P}}
\def\Z{{\mathbb Z}}
\def\R{{\mathbb R}} 
\def\ord{{\rm ord}}
\def\ssigma{\omega}
\def\f{\EuScript}
\def\e{\eqref}
\def\phi{{\varphi}}
\def\v{{\varepsilon}} 
\def\deg{{\rm deg\,}} 
\def\Det{{\rm Det}}
\def\dim{{\rm dim\,}} 
\def\Ker{{\rm Ker\,}} 
\def\Gal{{\rm Gal\,}}
\def\St{{\rm St\,}} 
\def\exp{{\rm exp\,}} 
\def\cos{{\rm cos\,}} 
\def\diag{{\rm diag\,}} 
\def\GCD{{\rm GCD }}
\def\LCM{{\rm LCM }}
\def\mod{{\rm mod\ }}
\def\Res{{\rm Res\ }}

\def\bp{\begin{proposition}}
\def\ep{\end{proposition}}
\newtheorem{zzz}{Theorem}
\newtheorem{yyy}{Corollary}
\def\bt{\begin{theorem}}
\def\et{\end{theorem}}
\def\be{\begin{equation}}
\def\bee{\begin{equation*}}
\def\la{\label}
\def\l{\lambda}
\def\m{\mu}
\def\ee{\end{equation}}
\def\eee{\end{equation*}}
\def\bl{\begin{lemma}}
\def\el{\end{lemma}}
\def\bc{\begin{corollary}}
\def\ec{\end{corollary}}
\def\pr{\noindent{\it Proof. }}
\def\note{\noindent{\bf Note. }}
\def\bd{\begin{definition}}
\def\ed{\end{definition}}
\def\e{\eqref}

\newtheorem{theorem}{Theorem}[section]
\newtheorem{lemma}{Lemma}[section]
\newtheorem{definition}{Definition}[section]
\newtheorem{corollary}{Corollary}[section]
\newtheorem{proposition}{Proposition}[section]

\begin{abstract}
In 1922 Ritt described polynomial solutions of the functional equation $P(f)=Q(g).$  
In this paper 
we describe solutions of the equation above in the case when $P,Q$ are polynomials while
$f,g$ are allowed to be arbitrary entire functions. In fact, we describe solutions 
of the more general functional equation $s=P(f)=Q(g),$ where $s,f,g$ are entire functions
and $P,Q$ are arbitrary rational functions. As an application we solve the problem of description 
of ``strong uniqueness polynomials'' for entire 
functions.
\end{abstract}

\maketitle

\section{Introduction} 
In this paper we describe all possible solutions of the functional equation  
\be \la{1} P\circ f=Q\circ g,\ee where $P,Q$ are polynomials, $f,g$ are entire functions, 
and the symbol $\circ$ denotes the superposition of functions, $f_1\circ f_2=f_1(f_2)$.
In fact we describe solutions 
of the more general functional equation
\be \la{11} s=P\circ f=Q\circ g,\ee where $s,f,g$ are entire functions
and $P,Q$ are arbitrary rational functions. 

Let us give several examples of solutions of \eqref{1}. First of all observe that 
for any polynomial $P$ and any entire function $f$ one can obtain a solution of \eqref{1}
setting \be \la{eq} Q=P\circ \alpha, \ \ \ \ g=\alpha^{-1}\circ f, \ee 
where $\alpha$ is a linear M\"obius transformation. Observe also that if $P,Q,f,g$ is a solution of \eqref{1}, then for any entire function $h$ and any polynomial $U$ the collection 
$$\hat P=U\circ P, \ \ \ \hat Q=U\circ Q, \ \ \ \hat f=f\circ h,  \ \ \ \hat g=g\circ h$$
also is a solution of \eqref{1}. 

In order to lighten the notation, in case if
rational functions $P,Q$ and entire functions $s,f,g$ such that  
\eqref{11} holds satisfy \eqref{eq} for some M\"obius transformation $\alpha$, we will say that that the decomposition $P\circ f$ of $s$ is equivalent to the decomposition $Q\circ g$. For equivalent decompositions we will use the notation $P\circ f\sim Q\circ g$. 

The simplest examples of solutions of \eqref{1} for which the decompositions $P\circ f$ and $Q\circ g$ are not equivalent are provided by 
polynomials. 
For example, we have $z^n\circ z^m=z^m\circ z^n$. 
More generally,
for any polynomial $R$ and $r\geq 0,$ $n\geq 1$ the equality \be \la{r1} z^n \circ z^rR(z^n)=z^rR^n(z) \circ z^n \ee holds.
Another examples of polynomial solutions of \eqref{1} are provided by the Chebyshev polynomials
$T_n$ defined by the equality \be \la{che} T_n(\cos z)=\cos nz.\ee Indeed, it follows from \eqref{che} that for any $m,n\geq 1$ we have: \be \la{r2} T_n\circ T_m=T_m\circ T_n.\ee 

The theory of functional decompositions of polynomials developed by Ritt \cite{ri} yields 
that actually any polynomial solution 
of \eqref{1} in a sense reduces either to \eqref{r1} or to \eqref{r2}. Namely, the following statement is true: if polynomials $P,Q,f,g$ satisfy \eqref{1}
then there exist polynomials $U$, $\tilde P,$ $\tilde Q,$ $\tilde f,$ $\tilde g,$ $h$ such that 
\be \la{rys}P=U \circ \tilde P, \ \ \ \ Q=U \circ \tilde Q, \ \ \ \ f=\tilde f \circ h, \ \ \ \ g=\tilde g \circ h, \ \ \ \ 
\tilde P\circ \tilde f=\tilde Q\circ \tilde g\ee 
and up to a possible replacement of $\tilde P$ by $\tilde Q$ and $\tilde f$ by $\tilde g$ either
\be \la{ru1} \tilde P\circ \tilde f\sim z^n \circ z^rR(z^n),  \ \ \ \ \ \ 
\tilde Q\circ \tilde g\sim  z^rR^n(z) \circ z^n,\ee 
where $R$ is a polynomial, $r\geq 0,$ $n\geq 1,$ and 
$\GCD(n,r)=1,$ or
\be \la{ru2} \tilde P\circ \tilde f\sim T_n \circ T_m, \ \ \ \ \ \ \tilde Q\circ \tilde g\sim T_m \circ T_n,\ee 
where $T_n,T_m$ are the corresponding Chebyshev polynomials with $n,m\geq 1,$ and $\GCD(n,m)=1.$

The simplest example of a solution of \eqref{1} with transcendental $f,g$ is provided by the equality 
$$\cos^2 z=1-\sin^2 z.$$ More generally, for any polynomial $S$ we have: 
\be \la{r3} z^2\circ \cos z\,S(\sin z)=(1-z^2)S^2(z)\circ \sin z.\ee 

The equality
\be \la{opi} T_n\circ \cos mz=T_m\circ \cos nz\ee
also is an example
of a solution of \eqref{1}. Nevertheless, in a sense this equality is a corollary of 
equality \eqref{r2}
since 
$$\cos mz =T_m \circ \cos z, \ \  \ \ \cos nz =T_n \circ \cos z.$$
On the other hand, for example the equality
$$-T_2\circ \cos \left(\frac{\pi}{2}+z\right )=T_2\circ \cos z$$ already can not be reduced in a similar way to \eqref{r1}, \eqref{r2}, or \eqref{r3}. 
More generally, for any $m,n\geq 1,$ $l>1$, and $0\leq k < nl$ we have: 
\be \la{r4} -T_{nl} \circ \cos\left( \frac{(2k+1)\pi}{nl}+mz\right)=T_{ml} \circ \cos(nz).\ee

Our first result states that up to one ``sporadic'' exception any solution of \eqref{1} can be reduced to \eqref{r1}, \eqref{r2}, \eqref{r3} or \eqref{r4}.

\pagebreak
\vskip 0.2cm 
\noindent{\bf Theorem A.} {\it Suppose that polynomials $P,Q$ and entire functions $f,g$ 
satisfy the equation $$P\circ f=Q\circ g.$$ Then 
there exist polynomials $F$, $\tilde P,$ $\tilde Q$ and entire functions 
$\tilde f,$ $\tilde g$, $t$ such that   
$$ P=F \circ \tilde P, \ \ \ \ Q=F \circ \tilde Q, \ \ \ \ 
f=\tilde f \circ t, \ \ \ \ g=\tilde g \circ t, \ \ \ \ \tilde P\circ \tilde f=\tilde Q\circ \tilde g$$
and, up to a possible replacement of $P$ by $Q$ and $f$ by $g$, 
one of the following conditions holds:
$$\tilde P\circ \tilde f\sim z^n \circ z^rR(z^n),  \ \ \ \ \ \ 
\tilde Q\circ \tilde g\sim  z^rR^n(z) \circ z^n,\leqno 1) $$ 
where $R$ is a polynomial, $r\geq 0,$ $n\geq 1,$  and 
$\GCD(n,r)=1;$ 
\vskip 0.01cm 
$$\tilde P\circ \tilde f\sim T_n \circ T_m, \ \ \ \ \ \ \tilde Q\circ \tilde g\sim T_m \circ T_n,\leqno 2)$$ where $T_n,T_m$ are the corresponding Chebyshev polynomials with $m,n\geq 1,$ and $\GCD(n,m)=1;$ 
\vskip 0.01cm
$$\tilde P\circ \tilde f\sim  z^2 \circ \,
\cos z\, S(\sin z),  \ \ \ \ \ \ \tilde Q\circ \tilde g\sim  (1-z^2)\, S^2(z)\circ \sin z,\leqno 3)$$
where $S$ is a polynomial;
\vskip 0.01cm
$$\tilde P\circ \tilde f\sim -T_{nl} \circ \cos\left( \frac{(2k+1)\pi}{nl}+mz\right), \ \ \ \ \ \ 
\tilde Q\circ \tilde g\sim T_{ml} \circ \cos(nz),\leqno 4)$$    
where $T_{nl},T_{ml}$ are the corresponding Chebyshev polynomials with $m,n\geq 1,$ \linebreak $\GCD(n,m)=1$, $l>1$, and $0\leq k < nl;$
\vskip -0.2cm 
$$\tilde P\circ \tilde f\sim (z^2-1)^3\circ  \left(
\frac{i\sin 2x+2\sqrt{2}\cos x}{\sqrt{3}}\right), \leqno 5) $$
$$\ \ \ \ \ \ \tilde Q\circ \tilde g\sim (3z^4-4z^3)\circ \left(\frac{i\sin 3x}{3\sqrt{2}}+\cos 2x 
+\frac{i\sin x}{\sqrt{2}}+\frac{2}{3}\right).$$
}
\vskip 0.2cm

Since a composition $P\circ f$ of a polynomial $P$ and an entire function $f$ is an entire function, 
the problem of description of solutions of 
equation \eqref{1} is a particular case of the problem of description of all possible ``double decompositions'' \eqref{11} of an entire function. 
Notice that different aspects of the theory of decompositions of entire functions were studied in many 
recent papers (see e. g. \cite{er}, \cite{l},
\cite{ny1}, \cite{ny2}, \cite{n2}, \cite{n}). However, 
this theory is still far from its completion.
In particular, there exist no results about double decompositions of entire functions similar to the results of Ritt.

Our next result 
describes solutions of equation \eqref{11} in case where the functions $P,Q$
are rational and 
at least one of them is 
not a polynomial. Together with Theorem A this
provides a complete description of solutions of equation \eqref{11} with rational $P,Q$ and entire 
$s,f,g.$

\pagebreak

\vskip 0.2cm 
\noindent{\bf Theorem B.} {\it Suppose that rational functions $P,Q$ and   
entire functions $s,f,g$ satisfy the equation 
$$s=P\circ f =Q\circ g.$$ Furthermore, suppose that at least one of the functions $P,Q$ is 
not a polynomial. 
Then there exist rational functions $F$, $\tilde P$, $\tilde Q$ and entire functions 
$\tilde f,$ $\tilde g,$ $t$ such that 
$$ P=F \circ \tilde P, \ \ \ \ Q=F \circ \tilde Q, \ \ \ \ 
f=\tilde f \circ t, \ \ \ \ g=\tilde g \circ t,\ \ \ \ \tilde P\circ \tilde f=\tilde Q\circ \tilde g$$
and, up to a possible replacement of $P$ by $Q$ and $f$ by $g$, 
one of the following conditions holds:
$$\tilde P\circ \tilde f\sim z^n \circ e^{rz}L(e^{nz}),  \ \ \ \ \ \ 
\tilde Q\circ \tilde g\sim  z^rL^n(z) \circ e^{nz},\leqno 1) $$ 
where $L$ is a Laurent polynomial, $r\geq 0,$ $n\geq 1,$ and 
$\GCD(n,r)=1;$ 
\vskip 0.01cm 
$$\tilde P\circ \tilde f\sim T_n \circ \cos(mz), \ \ \ \ \ \ \tilde Q\circ \tilde g\sim \frac{1}{2}\left(z^m+\frac{1}{z^m}\right) \circ e^{inz},\leqno 2)$$ where $T_n$ is the Chebyshev polynomial, $m,n\geq 1,$ and $\GCD(n,m)=1.$
}
\vskip 0.2cm

Yet another 
problem related to equation \eqref{1} is the problem of description of ``strong uniqueness polynomials'' for entire functions. Recall that a polynomial $P$ is called a strong uniqueness polynomial for entire functions if the equality 
\be \la{u} P\circ f=cP\circ g ,\ee where $f,g$ are entire functions and $c\in \C$, implies that $c=1$ and $f\equiv g.$
Such polynomials are closely related to the ``uniqueness range sets'' for entire functions 
and were studied in the recent papers \cite{cher},  \cite{fu1}, \cite{fu3}, \cite{sh}, \cite{yh} 
(see also the papers \cite{a0}, \cite{a}, \cite{fu2}, \cite{hy}, \cite{ly}, \cite{pak+} 
where the similar question was studied for meromorphic functions).

Usually, the problem of description of strong uniqueness polynomials for entire functions is studied under additional conditions of ``generic position"
imposed on the polynomial $P.$ Our last result based on Theorem A 
provides a complete solution of this problem in the general case.

\vskip 0.2cm 
\noindent{\bf Theorem C.} {\it A polynomial $P$ is not a strong uniqueness polynomial for entire functions if and only if
there exists a M\"obius transformation $\alpha$ such that 
either  
$$P=z^rR(z^n)\circ \alpha, \leqno 1)$$ where $R$ is a polynomial   
and $r\geq 0,$ $n>1,$ or 
$$P=F\circ T_l\circ \alpha,\leqno 2)$$ where $F$ is a polynomial, $T_l$ is  
the Chebyshev polynomial, and $l>1.$}
\vskip 0.2cm 

The paper is organized as follows. In the second section, using a result 
about
parametrizations of algebraic curves by entire functions obtained
in \cite{bn}, \cite{fs}, \cite{n} we relate the classification of solutions of equation \eqref{11} with
the classification of double decompositions of Laurent polynomials
into compositions of rational functions.
In the third section we review the papers \cite{bilu}, \cite{pak0}, \cite{pak}, \cite{z} where such a classification was obtained.
Finally, in 
the fourth section we prove Theorems A, B, and C.

\section{Reduction} Let $\f E:\,u(x,y)=0$ be an affine algebraic curve.
Recall, that a pair $f,g$ of functions meromorphic 
on a simply connected domain $D$ of $\C\P^1$ 
is called a meromorphic parametrization of $\f E$ 
on $D$ if for any point $z\in D$ which is not a pole of $f$ or $g$ the equality $u(f(z),g(z))=0$ holds and with finitely many exceptions any point
of $\f E$ is of the form $(f(z),g(z))$ for some $z\in D.$ 
If $D=\C$ and functions $f,g$ are entire then the corresponding parametrization is called 
entire.

Denote by $\hat {\f E}$ the desingularization of the curve $\f E.$ The general structure of meromorphic parametrizations of $\f E$ on $D$ is described by the following theorem (see \cite{bn}).

\bt \la{b} Let $f,g$ be a meromorphic parametrization of $\f E$ on 
$D$. Then there exist a holomorphic function $h:\, D\rightarrow \hat {\f E} ,$ meromorphic functions 
$U,V:\, \hat {\f E}  \rightarrow \C\P^1,$ and a finite set $S\subset \hat {\f E} $ 
such that $f=U\circ h,$ $g=V\circ h$ and the mapping
$(U,V):\, \hat {\f E}  \rightarrow \C\P^1\times \C\P^1$ is injective on $\hat {\f E} \setminus S.$ \qed
\et

The class of curves having a meromorphic 
parametrization on $\C$ is quite restrictive in view of the following classical Picard theorem \cite{pic}.

\bt If $\f E$ has a meromorphic parametrization on $\C$ then $\hat {\f E} $ has genus zero or one. \qed
\et  

Furthermore, for curves having an entire parametrization the following much more precise result holds (see \cite{fs}, \cite{bn}).

\bt \la{f} Let $f,g$ be an entire parametrization of $\f E$. Then $\hat {\f E} =\C\P^1$ and there exist an entire function $h$, rational functions $U,V$, and a finite set $S\subset \C\P^1$ 
such that $f=U\circ h,$ $g=V\circ h$ and the mapping 
$(U,V):\, \C\P^1\rightarrow \C\P^1\times \C\P^1$ is injective on $\f \C\P^1\setminus S$. \qed
\et

Theorem \ref{f} permits to relate the description of solutions of \eqref{11} with the description of solutions of the equation 
\be \la{mn} L=P\circ U=Q\circ V,\ee where $L$ is a Laurent polynomial and  
$P,Q,U,V$ are rational functions.

\bt \la{red} Suppose that rational functions $P,Q$ and entire functions $s,f,g$ satisfy equation \eqref{11}. 
Then 
there exist an entire function $h$, a Laurent polynomial $L$, and rational functions $U,V$ such that 
\be \la{re} s=L\circ h, \ \ \ f=U\circ h,\ \ \ g=V\circ h, \ \ \ L=P\circ U=Q\circ V.
\ee 
\et 

\pr If \eqref{11} holds, then $f,g$ is an entire parametrization of a factor of the 
algebraic curve
$$P_1(x)Q_2(y)-P_2(x)Q_1(y)=0,$$ 
where 
$P_1$, $P_2$ and $Q_1,$ $Q_2$ are pairs polynomials without common roots such that 
$$P=P_1/P_2, \ \ Q=Q_1/Q_2,$$
and hence by Theorem \ref{f} there exist an entire function $h$ and 
rational functions $U,V$ such that 
\be \la{y} f=U\circ h,  \ \ \ \ g=V\circ h.\ee Furthermore, it follows from  
$P\circ f=Q\circ g$ that 
$P\circ U=Q\circ V.$ Besides, clearly $s=L\circ h,$ where $L=P\circ U=Q\circ V.$

Since $s$ is an entire function, the equality $s=L\circ h$
implies that $h$ does not take any value in $\C$ which is a pole of $L$. On the other hand, by the Picard Little theorem an entire function may omit at most one value $a$ in $\C.$ Therefore,  
$L$ is either a polynomial or has at most one pole in $\C$ and in the last case this pole necessarily coincides with $a$. Therefore, replacing if necessary $h$ by $h-a$ and 
$L$ by $L\circ (z+a),$ without loss of generality we may assume that $L$ is a Laurent polynomial. \qed

\section{Decompositions of Laurent polynomials} 
In this section we review results concerning double decompositions of polynomials and Laurent polynomials. In accordance with the notation introduced above, if $P,Q,U,V$ are rational functions such that \eqref{mn} holds for some 
rational function $L$ and $$P=Q\circ \alpha,\ \ \ U= \alpha^{-1}\circ V$$ for some M\"obius transformation $\alpha$, then we will say that the decomposition $P\circ U$ of $L$ is 
equivalent to the decomposition $Q\circ V$ and will use for equivalent decompositions  
the notation $P\circ U\sim Q\circ V$.

The decomposition theory of polynomials was constructed by Ritt in the paper \cite{ri}. In particular, Ritt proved that if $\tilde P, \tilde Q, \tilde f, \tilde g$ are polynomials
satisfying $$\tilde P\circ \tilde f= \tilde Q\circ \tilde g$$ and 
such that \be \la{en} 
\GCD(\deg\tilde P,\deg\tilde Q)=1, \ \ \  \GCD(\deg\tilde f,\deg\tilde g)=1,\ee
then 
up to a possible replacement of $\tilde P$ by $\tilde Q$ and $\tilde f$ by $\tilde g$ either \eqref{ru1} or \eqref{ru2} holds. 

On the other hand, it was proved in \cite{en}, \cite{tor} (see also  \cite{mp}) 
that if $P,Q,f,g$ are arbitrary polynomials satisfying \eqref{1}, then there exist polynomials
$U, \tilde P, \tilde Q, \tilde f, \tilde g, h$ such that 
\be \la{en2} \deg U=\GCD(\deg P, \deg Q), \ \ \ \deg h=\GCD(\deg f, \deg g)\ee and equalities 
\eqref{rys} hold. Clearly, this result and the Ritt theorem 
taken together provide a full description of polynomial solutions of \eqref{1}. 

It is easy to see that the problem of description of polynomial solutions of equation \eqref{1}
essentially is equivalent to
the problem of description of algebraic curves of the form \be \la{cur} P(x)-Q(y)=0\ee having a factor of genus zero 
with one point at infinity. Indeed, 
if \eqref{cur} is such a curve, then the corresponding factor can be parametrized 
by some polynomials $f,g$ implying \eqref{1}, and vice versa if $P,Q,f,g$ is a polynomial solution of \eqref{1}, then \eqref{cur} 
has a factor of genus zero with one point at infinity.  
A more general problem
of description of curves \eqref{cur} having a factor
of genus $0$ with at most two points at infinity is closely related to the number theory
and was studied in the papers \cite{f1}, \cite{bilu}. In particular, in \cite{bilu}
an explicit list of such curves, defined over any field $k$ of characteristic zero, was obtained.

It was observed by the author several years ago that the problem of description of solutions of equation 
\eqref{mn}, where $L$ is a Laurent polynomial and $P,Q,U,V$ are rational functions, mostly reduces
to the above mentioned  problem of description of curves \eqref{cur} having a factor
of genus $0$ with at most two points at infinity.
Indeed, since a Laurent polynomial $L$ has at most two poles, it is easy to see that any decomposition of $L$ into a composition of two rational functions 
is equivalent either to a decomposition $A\circ L_1$, where $A$ is a polynomial and $L_1$ is a Laurent polynomial, or to a decomposition $L_2\circ B$, where $L_2$ is
a Laurent polynomial and $B=cz^d,$ $c\in \C,$ $d\geq 1.$
Therefore, a description of solutions of \eqref{mn}
reduces to a description of solutions of the following three equations:
\be \la{eg} A\circ L_1=B\circ L_2,
\ee where 
$A, B$ are polynomials and $L_1,L_2$ are Laurent polynomials, \be \la{2} A\circ L_1=L_2\circ z^d,\ee where 
$A$ is a polynomial and $L_1,L_2$ are Laurent polynomials, and
\be \la{3}L_1\circ z^{d_1}= L_2\circ z^{d_2},
\ee where $L_1,L_2$ are Laurent polynomials.

Equation \eqref{3} is very simple. Equation \eqref{2} is more complicated 
but still can be analysed quite easily in view of the presence of symmetries. 
Finally, if $A,B,L_1,L_2$ is a solution of equation \eqref{eg}, then 
the curve \be \la{cur1} A(x)-B(y)=0\ee
has a factor
of genus $0$ with at most two points at infinity and vice versa for any such a curve its factor of genus $0$ 
may be parametrized
by some Laurent polynomials.

A comprehensive self-contained theory of decompositions of Laurent polynomials was constructed in \cite{pak} where the equation
$$A\circ C=B\circ D$$ was studied 
in a more general setting involving holomorphic functions on compact Riemann surfaces.
In particular, in \cite{pak} were proposed new proofs of 
the Ritt theorem and the classification of curves \eqref{cur} having a factor of genus zero
with at most two points at infinity. Another approach to the problem may be found in \cite{z} where 
solutions of \eqref{eg} and \eqref{2} are deduced from 
results of \cite{bilu} and \cite{az}.

Below we collect necessary results from \cite{pak}. In order to lighten the notation set 
$$U_n(z)=\frac{1}{2} \left(z^n+\frac{1}{z^n}\right), \ \ \  V_n(z)=\frac{1}{2i} \left(z^n-\frac{1}{z^n}\right).$$
We start from the description of solutions of equation \eqref{eg}.

\bt \la{irrr} 
Suppose that polynomials $A,B$ and Laurent polynomials $L_1,$ $L_2$ satisfy the equation
$$A\circ L_1=B\circ L_2.$$
Then 
there exist polynomials $E,$ $\tilde A,$ $\tilde B$ and Laurent polynomials 
$W,$ $\tilde L_1,$ $\tilde L_2$ such 
that $$ A=E \circ \tilde A, \ \ \ \ B=E \circ \tilde B, \ \ \ \ 
L_1=\tilde L_1 \circ W, \ \ \ \ L_2=\tilde L_2 \circ W,\ \ \ \ \tilde A\circ \tilde L_1=\tilde B\circ \tilde L_2$$
and, up to a possible replacement of $A$ by $B$ and $L_1$ by $L_2$, one of the following conditions holds:
$$\tilde A\circ \tilde L_1\sim z^n \circ z^rR(z^n),  \ \ \ \ \ \ 
\tilde B\circ \tilde L_2\sim  z^rR^n(z) \circ z^n,\leqno 1) $$ 
where $R$ is a polynomial, $r\geq 0,$ $n\geq 1,$  and 
$\GCD(n,r)=1;$
\pagebreak
\vskip 0.01cm 
$$\tilde A\circ \tilde L_1\sim T_n \circ T_m, \ \ \ \ \ \ \tilde B\circ \tilde L_2\sim T_m \circ T_n,\leqno 2)$$ 
where $T_n,T_m$ are the corresponding Chebyshev polynomials with $m,n\geq 1,$ and $\GCD(n,m)=1;$ 
\vskip 0.01cm
$$\tilde A\circ \tilde L_1\sim  z^2 \circ U_1
S(V_1),  \ \ \ \ \ \ \tilde B\circ \tilde L_2\sim  (1-z^2)\,S^2 \circ V_1,\leqno 3)$$
where $S$ is a polynomial;
\vskip 0.01cm
$$\tilde A\circ \tilde L_1\sim -T_{nl} \circ U_m(\v z), \ \ \ \ \ \ 
\tilde B\circ \tilde L_2\sim T_{ml} \circ  U_n,\leqno 4)$$    
where $T_{nl},T_{ml}$ are the corresponding Chebyshev polynomials with $m,n\geq 1,$ $l>1$, $\varepsilon^{nlm}=-1$,
and $\GCD(n,m)=1;$ 
\vskip -0.2cm 
$$\tilde A\circ \tilde L_1\sim (z^2-1)^3\circ \left(
\frac{i}{\sqrt{3}}\,V_2+\frac{2\sqrt{2}}{\sqrt{3}}\,U_1\right), \leqno 5) $$
$$\ \ \ \ \ \ \tilde B\circ \tilde L_2\sim (3z^4-4z^3)\circ  \left(
\frac{i}{3\sqrt{2}}\,V_3+U_2
+\frac{i}{\sqrt{2}}\,V_1+\frac{2}{3} \right). \qed$$ 

\et

The expressions for $A$ and $B$ given in Theorem \ref{irrr} coincide with the ones given in \cite{bilu}, Theorem 9.3
(for $k=\C$) and \cite{pak}, Theorem 1.1 (see also Theorem 7.2 and Theorem 8.1 of \cite{pak}). The expressions for $L_1,L_2$ coincide with the ones given in \cite{pak}  
in 1), 2), 4) and slightly differ in 3) and 5). 
Say, in \cite{pak} in 3)  
we used for the curve $$x^2-(1-y^2)S(y)=0$$ 
the parametrization $$L_1= \frac{z^2-1}{z^2+1}
S(\frac{2z}{z^2+1}),  \ \ \ \ \ \ L_2=\frac{2z}{z^2+1}\,,$$ while now the parametrization
$$L_1=U_1S(V_1), \ \ \ \ \ L_2=V_1.$$ A similar change is made in 5).
Since the rational functions in the new parametrizations 
have the same degrees as the corresponding functions in the old parametrizations, the new parametrizations
can be obtained from the old ones by composing them with M\"obius transformations, 
and therefore such a change of parametrizations
does not affect the conclusion of the theorem.

The solutions of equation \eqref{2} (in the case where this equation does not reduce to \eqref{eg}) are described by the following theorem (see \cite{pak}, Theorem 6.4).

\bt \la{gop} Suppose that polynomials $A,$ $B$ and Laurent polynomials $L_1,$ $L_2$ (which are not polynomials)
satisfy the equation $$A\circ L_1=L_2\circ B.$$ 
Then there exist polynomials $E,$ $\tilde A,$ $\tilde B,$
$W$ and Laurent polynomials $\tilde L_1,$ $\tilde L_2$ such 
that 
$$A=E \circ \tilde A, \ \ \ \ L_2=E\circ \tilde L_2,  \ \ \ \ 
L_1=\tilde L_1 \circ W, \ \ \ \ B= \tilde B\circ W,\ \ \ \ \tilde A\circ \tilde L_1=\tilde L_2\circ \tilde B$$
and either 
\be \la{hy} \tilde A\circ \tilde L_1\sim z^{n}\circ z^rL(z^n),\ \ \ \ 
\tilde L_2\circ \tilde B \sim z^{r}L^{n}(z)\circ z^n,\ee
where $L$ is a Laurent polynomial, $r\geq 0,$ $n\geq 1,$ and $\GCD(r,n)=1,$ 
or 
\be \la{yh} 
\tilde A\circ \tilde L_1\sim  T_{n} \circ U_m,\ \ \ \ 
\tilde L_2\circ \tilde B \sim 
U_m\circ z^n,\ee
where $T_n$ is the Chebyshev polynomial, $n\geq 1,$ $m \geq 1,$ and $\GCD(m,n)=1.$ \qed
\et

Finally, solutions of equation \eqref{3} are described as follows
(see e.g. \cite{pak}, Lemma 6.1).

\bl\la{zc} 
Let $L_1,L_2$ be Laurent polynomials such that \eqref{3} holds for some $d_1,d_2\geq 1.$ 
Then there exists
a Laurent polynomial $N$ 
such that \be \la{ik} L_1=N\circ z^{D/d_1},\ \ \ \ \ \ L_2=N\circ z^{D/d_2},\ee
where $D=LCM(d_1,d_2).$ \qed
\el

\section{Proofs of Theorems A, B, and C}

\noindent{\it Proof of Theorem A.} Suppose that \eqref{1} holds for some entire functions $f,g$ and polynomials $P,Q$ and let $s$ be an entire function defined by equality \eqref{11}. By Theorem \ref{red} 
there exist an entire function $h$, a Laurent polynomial $L$, and rational functions $U,V$ such that equalities \eqref{re} hold. Furthermore, since
$P,Q$ are polynomials, the functions $U,V$ are Laurent polynomials. Therefore, setting
$$A=P, \ \ \ B=Q, \ \ \ L_1=U, \ \ \ L_2=V$$
we obtain \eqref{eg} and may apply Theorem \ref{irrr}. 

Observe that without loss of generality we may assume that the Laurent polynomial $W$ in Theorem \ref{irrr} equals $z.$
Indeed, if $W$ is a polynomial, then the function $W\circ h$ is clearly
entire, and we may simply replace $h$ by $W\circ h$.
On the other hand, if $W$ is 
not a polynomial, then 
$L$ is not a polynomial either, and it follows from \be \la{fty} s=L\circ h\ee that $h$ omits the value 0. Therefore, in this case the function $W\circ h$ is also entire and we may 
replace $h$ by $W\circ h$ as above.

If for the functions $A,B,L_1,L_2$ either
conclusion 1) or conclusion 2) of Theorem \ref{irrr} holds, then 
setting 
$$F=E, \ \ \tilde P=\tilde A, \ \ \tilde f=\tilde L_1, \ \ \tilde Q=\tilde B, \ \ \tilde g=\tilde L_2, \ \  t=h$$ 
we see that for $P,Q,f,g$ accordingly either 
conclusion 1) or conclusion 2) of Theorem A holds.

On the other hand, if for the functions $A,B,L_1,L_2$ one of conclusions 3), 4), 5) holds, then $L$ is not a polynomial. Therefore, 
$h$ omits the value 0 and hence there exists 
an entire function $w$ such that 
\be \la{ppo} h=e^{iz}\circ w.\ee 

If 3) holds, then it follows from \eqref{ppo}
taking into account the equalities  
\be \la{ega1} U_n\circ e^{iz}=\cos(nz), \ \ \ \  V_n\circ e^{iz}=\sin(nz),\ee that 
$$\tilde L_1\circ h=\cos zS(\sin z)\circ w,\ \ \ \tilde L_2\circ h=\sin z\circ w.$$ Therefore, 
setting $$F=E, \ \ \tilde P=\tilde A, \ \ \tilde f=\tilde L_1\circ e^{iz}, \ \ \tilde Q=\tilde B, \ \ \tilde g=\tilde L_1\circ e^{iz},\ \ \ t=w$$ we conclude that conclusion 3) of Theorem A holds.
Similarly, it follows from \eqref{ppo}, \eqref{ega1} that if for $A,B,L_1,L_2$ conclusion 5) of Theorem \ref{irrr} holds, then
for $P,Q,f,g$ conclusion 5) of Theorem A holds.

Finally, in order to prove that if for $A,B,L_1,L_2$ conclusion 4) of Theorem \ref{irrr} holds, then
for $P,Q,f,g$ conclusion 4) of Theorem A holds,
observe that the equality  
$\v^{nlm}=-1$ implies that $\v=e^{i\gamma},$ where 
$$\gamma=\frac{\pi }{nlm}(2r+1)$$ for some $r$, $0\leq r <nlm,$ and hence 
$$U_m\circ \v z \circ e^{iz}=U_m\circ e^{i(z+\gamma)}=
\cos\left( \frac{(2k+1)\pi}{nl}+mz\right),$$ 
for some $k,$ $0\leq k < nl.$ \qed
\vskip 0.2cm

\noindent{\it Proof of Theorem B.} It follows from Theorem \ref{red} 
that if \eqref{11} holds for some entire functions 
$s,f,g$ and rational functions $P,Q$, then  
there exist an entire function $h$, a Laurent polynomial $L$, and rational functions $U,V$ such that equalities \eqref{re} hold.
Furthermore, clearly
either \be \la{df} P\circ U\sim A\circ L_1, \ \ \ Q\circ V\sim L_2\circ z^d,\ee
where $A$ is a polynomial and $L_1,L_2$ are Laurent polynomials 
satisfying \eqref{2}, or
\be \la{df2} P\circ U\sim L_1\circ z^{d_1}, \ \ \ 
Q\circ V\sim L_2\circ z^{d_2},\ee 
where $L_1,L_2$ are Laurent polynomials satisfying \eqref{3}.

If \eqref{df} holds, then it follows from Theorem \ref{gop} that there exist rational functions $E,$ $\hat P,$ $\hat Q,$ $\hat U,$ $\hat V,$ 
$W$ such that 
\be \la{iio} P=E\circ \hat P,\ \  Q=E\circ \hat Q, \ \
U=\hat U \circ W, \ \ V=\hat V \circ W, \ \ 
\hat P\circ \hat U=\hat Q\circ \hat V\ee
and either 
\be \la{hy1} \hat P\circ \hat U\sim z^{n}\circ z^rL(z^n),\ \ \ \ 
\hat Q\circ \hat V \sim z^{r}L^{n}(z)\circ z^n,\ee
where $L$ is a Laurent polynomial, $r\geq 0,$ $n\geq 1,$ and $\GCD(r,n)=1,$ or 
\be \la{yh1} 
\hat P\circ \hat U\sim  T_{n} \circ U_m,\ \ \ \ 
\hat Q\circ \hat V\sim 
U_m\circ z^n,\ee
where $T_n$ is the Chebyshev polynomial,
$n\geq 1,$ $m \geq 1,$ and $\GCD(m,n)=1.$

Observe that since in both cases the function 
$\hat P\circ \hat U$ is a Laurent polynomial which is not a polynomial, it follows from $$L=E\circ \hat P\circ \hat U\circ W$$ 
and $L^{-1}\{\infty\}=\{0,\infty\}$ 
that 
$E$ is a polynomial and 
$W=cz^{\pm d},$ where $c\in \C^{\ast},$ $d\geq 1.$ 
Since \eqref{fty} implies that $h$ omits the value 0, replacing 
$h$ by $W\circ h$ we may assume that $W=z.$

Setting now accordingly to possibilities \eqref{hy1}, \eqref{yh1} 
either 
$$F=E, \ \ \tilde P=\hat P, \ \ \tilde f=\hat U\circ e^{z}, \ \ \tilde Q=\hat Q, \ \ \tilde g= \hat V\circ e^z,\ \ t=i w$$
or
$$F=E, \ \ \tilde P=\hat P, \ \ \tilde f=\hat U\circ e^{iz}, \ \ \tilde Q=\hat Q, \ \ \tilde g= \hat V\circ e^{iz},\ \ t=w,$$
where $w$ is an entire function such 
\eqref{ppo} holds,
we see that one of conclusions of Theorem B holds.

Similarly, if \eqref{df2} holds, then it is not hard to prove using Lemma \ref{zc} that there exist rational functions $\hat P,$ $\hat Q,$ $\hat U,$ $\hat V,$ a Laurent polynomial $E$, and a Laurent polynomial $W$ of the form $W=cz^{\pm d},$ $d=\GCD(d_1,d_2)$, $c\in \C^{\ast}$
such that \eqref{iio} holds and 
$$ \hat P\circ \hat U\sim z^{D/d_1}\circ z^{(d_1/d)z},\ \ \ \ 
\hat Q\circ \hat V \sim z^{D/d_2}\circ z^{(d_2/d)z}.$$
Furthermore, without loss of generality we may assume that $W=z$ and setting 
$$F=E, \ \ \tilde P=\hat P, \ \ \tilde f=\hat U\circ e^{z}, \ \ \tilde Q=\hat Q, \ \ \tilde g= \hat V\circ e^z,\ \ t=i w,$$ where $w$ is an entire function such \eqref{ppo} holds, 
we see that conclusion 1) of Theorem B holds with 
$L\equiv 1,$ $n=d_2/d,$ and $r=d_1/d.$ \qed
 
\vskip 0.2cm

\noindent{\it Proof of Theorem C.} Suppose that $P$ is a polynomial and $c$ is a complex number such that  
\be \la{hju} P\circ f=cP\circ g\ee
for some entire functions $f$ and $g$. Then for $P$, $f$, $g$, and $Q=cP$ one of conclusions of Theorem A holds. Furthermore, since $\deg P=\deg cP$ we have:
\be \la{ebs} \deg \tilde P=\deg \tilde Q.\ee

If 1) holds, then equality \eqref{ebs} together with the conditions $\GCD(r,n)=1$, $r\geq 0,$ $n\geq 1$ imply that  
\be \la{syr} P\circ f \sim cP\circ g.\ee 
The same is true if 2) holds. 

Further, it follows from \eqref{ebs} that 5) is impossible, while 4) may hold only if $m=n=1$. Finally, if 
3) holds, then necessarily $\deg S=0$, and it is easy to see that in this case 3) 
reduces to 4) with $l=2$, $m=n=1$, and $k=1.$ 
Summing up we see that equality \eqref{hju} implies that either $P\circ f$ is equivalent to $cP\circ g$ or \be \la{lkj} P=F\circ T_l\circ \alpha\ee for some polynomial $F,$ 
M\"obius transformation $\alpha$, and $l\geq 2.$

If \eqref{syr} holds, then there exist $a,b\in \C$ such that 
\be \la{ci} cP=P\circ (az+b).\ee Set $\hat P= P\circ (z-b)$. Then 
\be \la{ci1} c\hat P=\hat P \circ az\ee and 
the comparision of coefficients of polynomials in both parts of \eqref{ci1} 
implies that there exists an $n$th root of unity $\v$ such that $a=\v,$ $c=\v^r$, and  
\be \la{ert} \hat P=z^rR(z^n)\ee for some polynomial $R$ and $r\geq 0,$ $n>1$.
Therefore, \eqref{syr} implies that $$P=z^rR(z^n)\circ \alpha$$ for some M\"obius transformation 
$\alpha.$ Furthermore, any polynomial of such a from 
is not a uniqueness polynomial 
for entire functions since for any $n$th root of unity $\v$ distinct from 1 and any entire function $f$ we have: $$P\circ (\alpha^{-1}\circ f)=(\v^{n-r}P)\circ (\alpha^{-1}\circ \v f).$$

Finally, it is easy to see that any polynomial of form \eqref{lkj} also is not a strong uniqueness polynomial for entire functions
since for example for the functions
$$\alpha^{-1}\circ \cos\left(\frac{2\pi}{l}+z\right), \ \  \ g=\alpha^{-1}\circ  \cos z$$
the equality 
$$P\circ f=P\circ g$$ holds. \qed

\end{document}